\documentclass{amsart}

\usepackage{indentfirst}
\usepackage{amssymb}
\usepackage{amsmath}
\usepackage{graphics}
\usepackage{epsfig}
\usepackage[usenames]{color}
\usepackage{MnSymbol,wasysym}
\usepackage[utf8]{inputenc}
\usepackage{latexsym}
\usepackage{graphicx}
\usepackage[normalem]{ulem}
\usepackage{hyperref}
\usepackage[T1]{fontenc}
\usepackage{kantlipsum}
\usepackage[usenames,dvipsnames]{xcolor}
\usepackage[breakable, theorems, skins]{tcolorbox}
\usepackage{ mathrsfs }

\tcbset{enhanced}

\DeclareMathOperator*{\esssup}{ess\,sup}

\definecolor{light-gray}{gray}{0.95}

\catcode`@=11 \@addtoreset{equation}{section}
\renewcommand\theequation{\thesection.\@arabic\c@equation}
\catcode`@=12

\newcommand{\RR}{\mathbb{R}}

\expandafter\chardef\csname pre amssym.def
at\endcsname=\the\catcode`\@ \catcode`\@=11
\def\undefine#1{\let#1\undefined}
\def\newsymbol#1#2#3#4#5{\let\next@\relax
 \ifnum#2=\@ne\let\next@\msafam@\else
 \ifnum#2=\tw@\let\next@\msbfam@\fi\fi
 \mathchardef#1="#3\next@#4#5}
\def\mathhexbox@#1#2#3{\relax
 \ifmmode\mathpalette{}{\m@th\mathchar"#1#2#3}%
 \else\leavevmode\hbox{$\m@th\mathchar"#1#2#3$}\fi}
\def\hexnumber@#1{\ifcase#1 0\or 1\or 2\or 3\or 4\or 5\or 6\or 7\or 8\or
 9\or A\or B\or C\or D\or E\or F\fi}

\font\teneufm=eufm10 \font\seveneufm=eufm7 \font\fiveeufm=eufm5
\newfam\eufmfam
\textfont\eufmfam=\teneufm \scriptfont\eufmfam=\seveneufm
\scriptscriptfont\eufmfam=\fiveeufm

\catcode`\@=\csname pre amssym.def at\endcsname

\definecolor{light-gray}{gray}{0.95}

\newcommand{\eqn}{\begin{eqnarray}}
\newcommand{\een}{\end{eqnarray}}

\newtheorem {Theorem}  {Theorem}

\numberwithin{Theorem}{section}

\newtheorem{Remark}[Theorem]{Remark}

\newtheorem{Corollary}[Theorem]{Corollary}

\renewcommand{\div}{\mbox{div}}

 \newcommand{\ww}{\mbox{\boldmath $w$}}
 \newcommand{\vv}{\mbox{\boldmath $v$}}
\newcommand{\uu}{\mbox{\boldmath $u$}}
\newcommand{\zz}{\mbox{\boldmath $z$}}

\newcommand{\error}{\mbox{\boldmath $\varepsilon$}}

\begin{document}

\title[Strong alignment in $3D$ micropolar flows]{Strong alignment of micro-rotation and vorticity in $3D$ micropolar flows
}

\author[R.H. Guterres]{R.H. Guterres}
\address[R.H. Guterres]{Department of Pure and Applied Mathematics,  Universidade Federal do Rio Grande do Sul,  Porto Alegre, RS 91509-900, Brazil}
\email{rguterres.mat@gmail.com}

\author[W.G. Melo]{W.G. Melo}
\address[W.G.  Melo]{Department of Mathematics,   Universidade Federal de Sergipe,   S\~ao Crist\'ov\~ao, SE 49100-000, Brazil}
\email{wilberclay@gmail.com}

\author[C. J. Niche]{C.J.  Niche}
\address[C.J. Niche]{Departamento de Matem\'atica Aplicada, Instituto de Matem\'atica. Universidade Federal do Rio de Janeiro, CEP 21941-909, Rio de Janeiro - RJ, Brazil}
\email{cniche@im.ufrj.br}

\author[C. F.  Perusato]{C.F.  Perusato}
\address[C.F.  Perusato]{Department of Mathematics,  Universidade Federal de Pernambuco,  Recife, PE 50740-560, Brazil}
\email{cilon.perusato@ufpe.br}

\author[P.R.  Zingano]{P.R.  Zingano}
\address[P.R.  Zingano]{Department of Pure and Applied Mathematics,  Universidade Federal do Rio Grande do Sul,  Porto Alegre, RS 91509-900, Brazil}
\email{paulo.zingano@ufrgs.br}

\thanks{R. H. Guterres acknowledges support form FACEPE Grant BFP-0067-1.01/19.  W. Melo acknowledges support from Bolsa PQ CNPq - 309880/2021-1.  C.J.  Niche acknowledges support from PROEX CAPES - UFRJ.  C. F.  Perusato was partially supported by  CAPES PRINT-88881.311964/2018-01. }



\date{\today}

\begin{abstract} Rigid particles suspended on a micropolar fluid provide microstructure that coexists and interacts  with the local rotation of the fluid given by the vorticity.  In this work we prove that the particles' angular velocity and the vorticity strongly align for large times.  We provide average and supnorm estimates for the decay rate of the difference between these two vectors, which measures the alignment.
\end{abstract}

\maketitle

\begin{center}
{\it Dedicated in loving memory of Bella (2003 - 2022) and Sophie (2007 - 2018).}
\end{center}

\section{Introduction}

\subsection{Equations for micropolar fluids} Micropolar fluids are an important type of non-Newtonian fluids with velocity $\uu$ and pressure $p$  which have ``micro-structure'',  given by microscopic particles suspended or diluted in the fluid.  These rigid particles can rotate about their center of mass and have angular velocity $\ww$.  Since its introduction by Eringen \cite{MR0204005},  this fluid model has been extensely studied from a physical and mathematical point of view due to its importance in many applications.

We describe now the equations that model $3D$ micropolar fluids.  Given an initial incompressible velocity $\uu _0 \in L^{2}_{\sigma} $ and an initial micro-rotation $\ww_0 \in L^{2}$,  we consider a solution $(\uu,\ww)$,  in the Leray-Hopf  sense,  of the  micropolar fluid equations 

\begin{align}
\label{eqn:equation-u}
\partial_t \uu + \uu \cdot \nabla \uu + \nabla p & =  (\mu + \chi) \Delta \uu + 2 \chi \nabla \wedge \ww,  \\
\label{eqn:equation-w}
\partial_t \ww + \uu \cdot \nabla \ww & =  \nu \Delta \ww + \kappa \nabla (\nabla \cdot \ww) - 4 \chi \ww + 2 \chi \nabla \wedge \uu,  \\
\label{eqn:div-free}
\div \, \uu & =  0, 
\end{align}
see Eringen \cite{MR0204005}, Galdi and Rionero \cite{MR0467030},  Lukaszewicz \cite{MR1711268}, Yuan \cite{MR2419091}.  The coefficients  $\mu $ (kinematic viscosity),  $ \nu $ (angular viscosity)  and $ \chi $ (vortex or micro-rotation viscosity) are positive,  and $ \kappa $ (gyroviscosity) is nonnegative,  all of them assumed to be constant.   Leray-Hopf (or simply Leray) solutions are global mappings

\begin{displaymath}
\zz = (\uu,\ww) \in C_{\sf w} \left([0, \infty), L^{2}_{\sigma} \times L^{2}  \right) \cap L^{2} \left( (0, \infty), 
\dot{H}^{1} \times \dot{H}^{1} \right)
\end{displaymath}
with $\zz_0 (x) = \zz (x,0) = (\uu,\ww) (x,0) = (\uu_0 (x), \ww_0 (x))$  that satisfy the equations in weak sense for $ t > 0$ and also the energy estimate

\begin{equation}
\label{eqn:energy-inequality}
\Vert \zz (t) \Vert _{L^2} ^2 + \int _s ^t \left( \mu \Vert D \uu (\tau) \Vert _{L^2} ^2 + \nu \Vert D \ww (\tau) \Vert _{L^2} ^2\right) \, d \tau \leq \Vert \zz (s) \Vert _{L^2} ^2
\end{equation}
for all $t > s$ and a.e. $s > 0$ (see Section \ref{notations} for notation).    The existence of Leray solutions for the equations \eqref{eqn:equation-u} - \eqref{eqn:div-free} and other similar dissipative systems is well known, 
but their uniqueness and exact regularity properties are still open except for small initial data in suitable spaces, see Boldrini and Rojas-Medar \cite{MR1666509},  Galdi and Rionero \cite{MR0467030},   Lukaszewicz \cite{MR1711268},  Temam \cite{MR769654},  Yuan \cite{MR2419091}.  This is the case, for example, of small initial data in $ H_{\sigma}^{1} \times H^{1} $, as  Serrin-like regularity conditions (see Dong and Chen \cite{MR2572698}, Yuan \cite{MR2778615})  and standard calculations lead to global smooth solutions $\zz \in C^{\infty}(\mathbb{R}^{3} \times (0, \infty))$ such that $\zz \in  C \left( (0, \infty),  H^{m}  \times H^{m} \right)$ for all $ m \geq 0 $,  provided

\begin{equation}
\label{eqn:smallness-condition}
\Vert (\uu_0, \ww_0)  \Vert ^{\frac{1}{2}} _{L^2} \Vert (D \uu_0, D \ww_0)  \Vert ^{\frac{1}{2}} _{L^2}  \leq 3.182 \cdot \min \, \{\mu, \nu \},
\end{equation}
see \eqref{eqn:smallness}.  From this and \eqref{eqn:energy-inequality} it follows that for any $\zz_0 = (\uu_0, \ww_0) \in L_{\sigma}^{2} \times L^{2}$,  there is a $t_{\ast} \geq 0$ such that  $\zz \in C^{\infty}\left( \mathbb{R}^{3} \times \left( t_{\ast}, \infty \right) \right)$ and $\zz \in  C \left( (t_{\ast}, \infty),  H^{m}  \times H^{m} \right)$.  Moreover, there exists $t_{\ast \ast} \geq t_{\ast}$ such that 

\begin{equation}
\label{eqn:monotonicity-time}
0 \leq\ t_{\ast} \leq t_{\ast\ast} \leq K \cdot \left( \min \{\mu, \nu \} \right)^{- 5} \, \Vert \zz_0 \Vert ^4 _{L^2}, \qquad K < 0.005,
\end{equation}
for which $ \Vert D \zz (t) \Vert _{L^2}$ is monotonically decreasing in $(t_{\ast \ast}, \infty)$,  see Theorem \ref{monotonicity} in Section \ref{review}.  From \eqref{eqn:energy-inequality} and \eqref{eqn:monotonicity-time} we immediately obtain 

\begin{displaymath}
\lim_{t\,\rightarrow\,\infty}  t^{\frac{1}{2}} \, \Vert D \zz (t) \Vert _{L^2} = 0. 
\end{displaymath}
From this it follows that for arbitrary $\zz_0 \in L_{\sigma}^{2} \times L^{2}$ we have that 

\begin{equation}
\label{eqn:decay-dmu-not-sharp}
\lim_{t\,\rightarrow\,\infty}  t^{\frac{m}{2}} \, \Vert D ^m \uu  (t) \Vert _{L^2} = 0
\end{equation}
and 

\begin{equation}
\label{eqn:decay-dmw-not-sharp}
\lim_{t\,\rightarrow\,\infty}  t^{\frac{m+1}{2}} \, \Vert D ^m \ww  (t) \Vert _{L^2} = 0
\end{equation}
for every $ m \geq 0 $, see Guterres, Melo, Nunes and Perusato \cite{MR3906315} and Guterres, Nunes and Perusato \cite{MR3854334}.  The extra $t^{\frac{1}{2}}$ factor in \eqref{eqn:decay-dmw-not-sharp} with respect to \eqref{eqn:decay-dmu-not-sharp} seems to have been first obtained in the incompressible case in Guterres, Melo, Nunes and Perusato \cite{MR3906315},  Guterres, Nunes and Perusato \cite{MR3854334} (for compressible flows,  see Liu and Zhang \cite{MR3473451}),  as a consequence of the damping term $4 \chi \ww$ in the equation \eqref{eqn:equation-w}.  When the initial data are also in $ L^1$, the corresponding results were obtained by Cruz, Perusato,  Rojas-Medar and P. Zingano \cite{MR4358138}  by combining the monotonicity approach of Guterres, Niche,  Perusato and P. Zingano \cite{MR4546674} and Kreiss, Hagstrom, Lorenz and P. Zingano \cite{MR1994780} with the Fourier Splitting Method of M.E. Schonbek \cite{MR775190}, \cite{MR837929}.  Using the  Decay Character Theory of Bjorland and M.E. Schonbek \cite{MR2493562} upper and  lower bounds for $\Vert \zz (t) \Vert _{L^2}$ were obtained by  Niche and Perusato \cite{MR4379088} for some classes of initial data.   However,  an individual sharper lower bound for $\Vert \ww(t) \Vert _{L^2}$ was not obtained there.  The linear damping term gives an exponential  contribution  to the decay,  which prevents obtaining that lower estimate.  This has now become possible in view of the  strong decay properties of the ``error term'' $\error (x,t)$ given in \eqref{eqn:error} below.  It will follow from our results   that the speed-up gain by $ t^{\frac{1}{2}}$  observed in \eqref{eqn:decay-dmw-not-sharp} is indeed optimal and is rooted in an intimate relation between micro-rotation and flow vorticity that is an important effect of vortex viscosity.

Gay-Balmaz, Ratiu and Tronci  \cite{MR3116004} established the equivalence of Eringen's micropolar theory  and the  Ericksen - Leslie director theory for the dynamics of nematic liquid crystals developed by Ericksen \cite{MR137403}, \cite{MR1079183} and Leslie \cite{Leslie1979TheoryOF}.   Therefore, the results provided here might be  extended to other contexts.

\subsection{Physical motivation} We present now  the ideas that motivate our study.  From \eqref{eqn:equation-w},  it is clear that the fluid advects its microstructure through the term $\uu \cdot \nabla \ww$.  Note that this is the only term in \eqref{eqn:equation-u} and \eqref{eqn:equation-w} in which $\uu$ and $\ww$ interact directly.  As any vector field $\uu$ can be infinitesimally decomposed as a translation,  a deformation through shear and a rotation around an axis  given by the vorticity $\frac{1}{2} \nabla \wedge \uu$ (see pages 19  - 21,  Chorin and Marsden \cite{MR1218879}), we would like to understand how the vorticity and the micro-rotation compare to each other at a common ``micro'' scale.  This makes sense from the dimensional point of view,  as both  vectors  are proportional to the vortex viscosity $\chi > 0$ in  \eqref{eqn:equation-w}.  From Braz e Silva,  Cruz,  Freitas and P. Zingano \cite{MR3955606},  Cruz and Novaes \cite{CruzNovais},  Niche and Perusato \cite{MR4379088},  Cruz, Perusato, Rojas-Medar and P.  Zingano \cite{MR4358138},   we know that 

\begin{displaymath}
\Vert \ww(t) \Vert _{L^2}, \Vert \nabla \wedge \uu (t) \Vert _{L^2} = O \left( t^{-\left( \alpha + \frac{1}{2} \right)} \right),  
\end{displaymath}
for the same $\alpha > 0$,  which depends on the regularity of the initial data.   Moreover,  for a certain  range of values of $\alpha$,  we also have

\begin{displaymath}
\Vert \ww(t) \Vert _{L^2}, \Vert \nabla \wedge \uu (t) \Vert _{L^2} \geq C \, t^{-\left( \alpha + \frac{1}{2} \right)} 
\end{displaymath}
see Niche and Perusato \cite{MR4379088}.  Then,  in order to be  able to provide a more precise meaning to the qualitative assertion $ \ww \approx \frac{1}{2} \nabla \wedge \uu$,  we need to find a good estimate for the difference

\begin{equation}
\label{eqn:error}
\error (x,t) = \ww(x,t) -  \frac{1}{2} \nabla \wedge \uu (x,t).
\end{equation}

\subsection{Main results} \label{main-results} We now describe our main results. We are concerned throughout this paper with Leray solutions $\mbox{\boldmath $z$}(x,t) = \left(\uu (x,t),  \ww(x,t) \right)$ such that the velocity $\uu$ satisfies

\begin{equation}
\label{eqn:upper-bound-u}
\Vert \uu(t) \Vert _{L^2} \leq C_0 \, t ^{- \alpha},  \qquad \forall \, t \geq T_0,
\end{equation}
for some constants $\alpha, T_0,  C_0 = C_0 \left( \Vert \mbox{\boldmath $z$}_0 \Vert _{L^2}, \alpha, T_0 \right) \geq 0 $.   For the existence of solutions satisfying \eqref{eqn:upper-bound-u} with $ \alpha > 0 $,  see Braz e Silva,  Cruz,  Freitas and P. Zingano \cite{MR3955606},  Braz e Silva, Guterres, Perusato and P. Zingano \cite{silva2022high},  Dong and Chen \cite{MR2338669},  Niche and M.E. Schonbek \cite{MR3355116}.

We first state our results concerning upper bounds for $\uu,  \ww,  \error$ and its derivatives.

\begin{Theorem}
\label{decay-du-dw-upper}
If \eqref{eqn:upper-bound-u} is valid for some $\alpha \geq 0$,  then we have that for $m \geq 1$

\begin{equation}
\label{eqn:upper-bound-du}
\Vert D^m \uu(t) \Vert _{L^2}  \leq C \, \mu ^{- \frac{m}{2}} \, t ^{- \left( \alpha + \frac{m}{2} \right)},  \qquad \forall \, t \geq T_m,
\end{equation}
where $C = C (C_0, m, \alpha)$ and $T_m >  t_{\ast \ast}$ sufficiently large,   see Remark \ref{dependence-times-constants} for the dependence of $T_m$ on the constants appearing in the problem.  Moreover, we have that for every $m \geq 0$

\begin{equation}
\label{eqn:upper-bound-dw}
\Vert D^m \ww(t) \Vert _{L^2} \leq  C \, \mu ^{- \frac{m + 1}{2}} \,  t ^{- \left( \alpha + \frac{m + 1}{2} \right)},  \qquad \forall \, t \geq T'_m,
\end{equation}
where $C = C (C_0, m+1, \alpha)$ and  $T'_m >  t_{\ast \ast}$ sufficiently large,  see Remark \ref{dependence-times-constants}. 
\end{Theorem}

\begin{Theorem}
\label{decay-error-upper}
If \eqref{eqn:upper-bound-u} is valid for some $\alpha \geq 0$,  then we have that for $m \geq 0$

\begin{equation}
\label{eqn:upper-bound-error}
\Vert D^m \error (t) \Vert _{L^2} \leq C_1 \, |\mu - \nu| \, t ^{- \left( \alpha + \frac{m + 3}{2}\right)} + C_2 \, t^{- \left( 2 \alpha + \frac{m+3}{2} + \frac{1}{4} \right)}, \qquad \forall \, t \geq T''_m, 
\end{equation}
where $C_i$ are $C_i = C_i (C_0, m, \alpha, \mu, \chi), i = 1,2$ and  $T''_m >  t_{\ast \ast}$ sufficiently large,  see Remark \ref{dependence-times-constants}. 
\end{Theorem}

\begin{Remark} \label{constant-mu}
We emphasize the fact that the constants in \eqref{eqn:upper-bound-du} and \eqref{eqn:upper-bound-dw} depend on $\mu$,  but not on $\chi$  and $\nu$,  see  also Section \ref{sketch}.  
\end{Remark}

\begin{Remark} \label{dependence-times-constants} In Theorems \ref{decay-du-dw-upper} and \ref{decay-error-upper},  the times $T_m, T_m ', T_m '' > t_{\ast \ast}$ for which the estimates are valid,  depend  on the constants $ C_0, \alpha,$ and  $T_0$ in \eqref{eqn:upper-bound-u}, on the order $m$ of derivation,  and on the constants $\mu, \chi, \kappa, \nu$ in \eqref{eqn:equation-u} and \eqref{eqn:equation-w}. 
\end{Remark}

\begin{Remark}
From \eqref{eqn:error}  and Theorem \ref{decay-error-upper} we obtain that for $\mu \neq \nu$ and $m \geq 0$

\begin{displaymath}
\Vert D^m \left ( \nabla \cdot \ww \right) (t) \Vert _{L^2} = O \left( t^{- \left( \alpha + 2 + \frac{m}{2} \right)}\right), 
\end{displaymath}
while for $\mu = \nu$ and $m \geq 0$

\begin{displaymath}
\Vert D^m \left ( \nabla \cdot \ww \right) (t) \Vert _{L^2} = O \left( t^{- \left( 2\alpha + \frac{9}{4} + \frac{m}{2} \right)}\right). 
\end{displaymath}
From Theorem \ref{decay-du-dw-upper} we obtain for $m \geq 0$

\begin{displaymath}
\Vert D^m \left ( \nabla \wedge \ww \right) (t) \Vert _{L^2} = O \left( t^{- \left( \alpha + 1 + \frac{m}{2} \right)}\right). 
\end{displaymath}
\end{Remark}

Estimates \eqref{eqn:upper-bound-du}, \eqref{eqn:upper-bound-dw} and \eqref{eqn:upper-bound-error} seem to imply that the decay of $\error$ may be significantly faster than that of $\nabla \wedge \uu$ and $\ww$.  However,   \eqref{eqn:upper-bound-du} and \eqref{eqn:upper-bound-dw} are only upper bounds and do not preclude the possibility of faster decay.   The assumption

\begin{equation}
\label{eqn:lower-bound-u}
\Vert \uu(t) \Vert _{L^2} \geq \widetilde{C}_0 \, t^{- \alpha},  \qquad \forall \, t > t_0,
\end{equation} 
for some constants $\alpha, t_0,  \widetilde{C}_0 = \widetilde{C}_0 \left( \Vert \mbox{\boldmath $z$}_0 \Vert _{L^2}, \alpha, t_0 \right) \geq 0 $,  will lead to lower bounds analogous to  \eqref{eqn:upper-bound-du} and  \eqref{eqn:upper-bound-dw}.     Existence of solutions for which \eqref{eqn:upper-bound-u} and \eqref{eqn:lower-bound-u} hold was proved by Dong and Chen \cite{MR2338669} and Niche and Perusato \cite{MR4379088}, for the case $0 < \alpha < \frac{1}{2}$.

\begin{Remark} Assumptions  \eqref{eqn:upper-bound-u} and \eqref{eqn:lower-bound-u} seem to be natural,  as recently Brandolese,  Perusato and P. Zingano \cite{https://doi.org/10.1112/blms.12912} proved that solutions to the Navier-Stokes equations in $\RR ^n$ for which $\Vert \uu (t) \Vert _{L^2} = O \left( t ^{- \alpha } \right)$,  with $0 < \alpha < \frac{n+2}{4}$,  generically have the same lower bound.  
\end{Remark} 

We now state our results on lower bounds for the derivatives of $\uu$ and discuss their consequences.

\begin{Theorem}
\label{lower-bound-du}
If \eqref{eqn:upper-bound-u} and \eqref{eqn:lower-bound-u} hold for some $\alpha > 0$, then we have that  for $m \geq 1$

\begin{equation}
\label{eqn:lower-bound-du}
\Vert D^m \uu (t) \Vert _{L^2} \geq C \, \mu ^{- \frac{m}{2}} \, t^{- \left( \alpha + \frac{m}{2} \right)},  \qquad \forall \, t > t_m,
\end{equation}
where $C = C(m, \alpha, C_0, \widetilde{C}_0)$ and $t_m > t_{\ast \ast}$  sufficiently large,   see Remark \ref{dependence-times-constants}. 

\end{Theorem}

\begin{Remark}
The time $t_m$,  analogously to the ones described in Remark \ref{dependence-times-constants}, depends on the constants listed there,  but also on $\widetilde{C}_0$ and $t_0$.
\end{Remark}

\begin{Remark}
As $\Vert D^m  \nabla \wedge \uu (t)  \Vert _{L^2} = \Vert D^{m+1} \uu (t) \Vert _{L^2}$,  from \eqref{eqn:error},  \eqref{eqn:upper-bound-error} and \eqref{eqn:lower-bound-du}  we obtain that for $m \geq 0$

\begin{equation}
\label{eqn:lower-bound-dw}
\Vert D^m \ww (t) \Vert _{L^2} \geq C \, \mu^{- \frac{m+1}{2}} \, t^{- \left( \alpha + \frac{m+1}{2} \right)} + O \left( t^{- \left(\alpha + \frac{m+3}{2} \right) } \right),  \qquad \forall \, t \gg 1, 
\end{equation}
hence \eqref{eqn:upper-bound-u} and \eqref{eqn:lower-bound-u} yield upper and lower bounds for $\ww$ as well. 
\end{Remark}

\begin{Remark}
From \eqref{eqn:upper-bound-du}, \eqref{eqn:upper-bound-dw},  \ref{eqn:lower-bound-du} and \ref{eqn:lower-bound-dw}   and the Gagliardo-Nirenberg inequalities

\begin{displaymath}
\Vert D^m f \Vert _{L^{\infty}} \leq C \, \Vert D^m f \Vert _{L^2}  ^{\frac{1}{4}} \, \Vert D^{m+2} f \Vert _{L^2}  ^{\frac{3}{4}}, \quad \Vert D^m f \Vert _{L^{\infty}} \geq C \, \Vert D^{m+2} f \Vert _{L^2}  ^{\frac{3}{2}}  \, \Vert D^{m+3} f \Vert _{L^2}  ^{- \frac{1}{2}}
\end{displaymath}
we obtain the $L^{\infty}$ estimates

\begin{align*}
C_1  \, \mu ^{- \left( \frac{m}{2} + \frac{3}{4} \right)} \, t^{- \left( \alpha + \frac{m}{2} + \frac{3}{4} \right)} & \leq \Vert D^m \uu (t) \Vert _{L^{\infty}} \leq C_2 \, \mu ^{- \left( \frac{m}{2} + \frac{3}{4} \right)} \, t^{- \left( \alpha + \frac{m}{2} + \frac{3}{4} \right)},  \\ C_1  \, \mu ^{- \left( \frac{m+1}{2} + \frac{3}{4} \right)} \, t^{- \left( \alpha + \frac{m+1}{2} + \frac{3}{4} \right)} & \leq \Vert D^m \ww (t) \Vert _{L^{\infty}} \leq C_2  \, \mu ^{- \left( \frac{m+1}{2} + \frac{3}{4} \right)} \, t^{- \left( \alpha + \frac{m+1}{2} + \frac{3}{4} \right)}.
\end{align*}
\end{Remark}

Now,  we consider the triangle with sides $\ww,  \frac{1}{2} \nabla \wedge \uu$ and $\error = \ww - \frac{1}{2} \nabla \wedge \uu$.   As from  \eqref{eqn:upper-bound-du}, \eqref{eqn:upper-bound-dw},   \eqref{eqn:lower-bound-du} and \eqref{eqn:lower-bound-dw} we have that for large enough $t > t_{\ast \ast}$

\begin{displaymath}
C_1 \, t^{- \left(\alpha + \frac{1}{2}  \right)} \leq \Vert \nabla \wedge \uu (t)  \Vert _{L^2},   \Vert \ww (t) \Vert _{L^2} \leq  C_2 \, t^{- \left(\alpha + \frac{1}{2}  \right)},
\end{displaymath}
and from \eqref{eqn:upper-bound-error}

\begin{displaymath}
\Vert \error (t) \Vert _{L^2} \leq C \, t ^{- \left(\alpha + \frac{3}{2}  \right)}, 
\end{displaymath}
then $\ww$ and  $ \frac{1}{2} \nabla \wedge \uu$ tend to become aligned in $L^2$,  i.e.  on average they are colinear and have the same direction.  Moreover,  also from \eqref{eqn:upper-bound-du}, \eqref{eqn:upper-bound-dw},   \eqref{eqn:lower-bound-du} and \eqref{eqn:lower-bound-dw}

\begin{displaymath}
C_1 \, t^{- \left(\alpha + \frac{m+1}{2}  \right)} \leq \Vert D^m  \nabla \wedge \uu (t)  \Vert _{L^2},   \Vert D^m \ww (t) \Vert _{L^2} \leq  C_2 \, t^{- \left(\alpha + \frac{m+1}{2}  \right)},
\end{displaymath}
and from \eqref{eqn:upper-bound-error}

\begin{displaymath}
\Vert D^m \error (t) \Vert _{L^2} \leq C \, t ^{- \left(\alpha + \frac{m + 3}{2}  \right)}, 
\end{displaymath}
hence their derivatives are also $L^2$ aligned.  Using these estimates and Gagliardo-Nirenberg inequalities  we prove that this strong alignment also holds in the $L^{\infty}$ norm.

\begin{Corollary}
Suppose \eqref{eqn:upper-bound-u} and \eqref{eqn:lower-bound-u} hold for some $\alpha > 0$.  Then,  if  $\mu \neq \nu$,  for all $m \geq 0$

\begin{displaymath}
\Vert D^m \error (t) \Vert _{L^{\infty}} \leq C \, |\mu - \nu|  t^{- \left( \alpha + \frac {9}{4} + \frac{m}{2} \right)},
\end{displaymath}
and if $\mu = \nu$, then

\begin{displaymath}
\Vert D^m \error (t) \Vert _{L^{\infty}} \leq  C \, t^{-  \left( 2 \alpha + \frac{5}{2} + \frac{m}{2}   \right)}.
\end{displaymath}
\end{Corollary}

\subsection{Sketch of the proof} \label{sketch}We describe now the main ideas behind the proof of our estimates.  As solutions are eventually regular,  i.e.  they are $C^{\infty}$ smooth for large enough times,  all of our estimates are rigorous.  One of the main tools used in our proofs is the monotonicity property for $\Vert D \zz (t) \Vert _{L^2}$ in Theorem \ref{monotonicity},  which can be extended to  $\Vert D^m \zz (t) \Vert _{L^2}$,  see Theorem \ref{monotonicity-dm}.  From these,  the decay estimate  

\begin{displaymath}
\Vert D^m \zz (t) \Vert _{L^2} \leq C \, t^{- \frac{m}{2}},  m \geq 1, 
\end{displaymath} 
can be proved,  see Theorems \ref{decay-dz} and \ref{decay-dmz}.  This,  in turn,  leads to

\begin{displaymath}
\Vert D ^m \ww (t) \Vert _{L^2} \leq C \,  t ^{-\frac{m+1}{2}},  m \geq 0, 
\end{displaymath} 
as seen in Theorem \ref{decay-dmw}.  Note that the constants above depende on $m, C_0,  \mu, \nu, \chi$.

The proof of these preliminary estimates leads us to notice that faster decay rates may be obtained if we have further  information concerning  $\Vert \zz (t) \Vert _{L^2}$. This is due to  the finer estimate

\begin{equation}
\label{eqn:two-estimates-depending-on-z}
\Vert D^m \zz (t) \Vert _{L^2} \leq C \, \Vert \zz (t/2) \Vert _{L^2} \, t ^{- \frac{m}{2}},  \quad \Vert D^m \ww (t) \Vert _{L^2} \leq \, C \Vert \zz (t/2) \Vert _{L^2} \, t ^{- \frac{(m + 1)}{2}},
\end{equation}
see Remark \ref{bootstrap-dmz-dmw}.  From assumption \eqref{eqn:upper-bound-u},  i.e.  $\Vert \uu (t) \Vert_{L^2}  \leq C_0 \, t^{- \alpha}, \alpha > 0$,  we deduce that $\Vert \ww (t) \Vert _{L^2} \leq C \, t^{- \alpha}$,  which leads to Theorem \ref{decay-du-dw-upper-gamma},  that is, for $m \geq 0$

\begin{displaymath}
\Vert D^m \uu (t) \Vert _{L^2} \leq C \, t^{- \left( \alpha + \frac{m}{2} \right)}, \quad \Vert D ^m \ww (t) \Vert _{L^2} \leq C \,  t ^{- \left( \alpha +  \frac{m+1}{2} \right)}.
\end{displaymath}
This is a first version of Theorem \ref{decay-du-dw-upper},  but now the constant $C$ depends on $\gamma = \min \{ \mu, \nu \}$, instead of just $\mu$. 

A second consequence of Remark \ref{bootstrap-dmz-dmw} and \eqref{eqn:two-estimates-depending-on-z} is Theorem \ref{estimate-error-first},  which provides a preliminary estimate for the decay of $\error$, 

\begin{displaymath}
\Vert D^m \error (t) \Vert _{L^2} \leq C t ^{- \left( \alpha + \frac{m+3}{2} \right)}.
\end{displaymath}
This result,  together with a new monotonicity estimate \eqref{eqn:zeta-m} in Theorem  \ref{monotonicity-zeta},  leads to the proof of Theorem \ref{decay-du-dw-upper}.  In turn,  the estimates in this last result and  in Theorem \ref{estimate-error-first},  lead to the proof of Theorem \ref{decay-error-upper}.  Finally,  \eqref{eqn:zeta-m} in Theorem  \ref{monotonicity-zeta}, Theorem \ref{decay-du-dw-upper} and   Theorem \ref{estimate-error-first}, imply Theorem \ref{lower-bound-du}.  

\subsection{Organization of this paper} This article is organized as follows.  In Section 2, we gather notation, definitions and results needed for our proofs.  Although this material is basically known,  some proofs are new and improvements are offered.  In Section \ref{proofs} we prove our main results. 

\section{Mathematical preliminaries}
\label{preliminaries}

\subsection{Summary of notations}
\label{notations}

Throughout this text $\langle \vv ,  \ww \rangle $ indicates the standard inner product of two vectors $\vv$ and $\ww$ in $\RR^3$,  and $ | \cdot | $ is used to denote the absolute value of a scalar or the Euclidean norm of a vector.   The usual Lebesgue spaces are  denoted by $ L^{p}, 1 \leq p \leq \infty $; $\vv = (v_{1}, v_{2}, v_{3})  \in L^{p} $ means that $ v_{i} \in L^{p} $ for every $ 1 \leq i \leq 3 $.  By $ \| \cdot \|_{L^{p}}$ we mean  the  usual $L^{p}$ norm; for vector functions, 

\begin{displaymath}
\Vert \vv \Vert ^p _{L^p} = \Vert (v_{1},  \cdots,  v_{k}) \Vert ^p _{L^p} = \sum_{i=1} ^{3} \int _{\RR ^3} | v_{i}(x) | ^{p} \, dx
\end{displaymath}
if $1 \leq p < \infty $,  and
\begin{displaymath}
\Vert \vv \Vert _{L^{\infty}} = \esssup _{x \in \RR^3} \{ |\left(v_1 (x), \cdots, v_3 (x) \right| \},
\end{displaymath}
if $ p = \infty$.   Similarly,  if $1 \leq p < \infty$  and $v = (v_1, v_2,  v_3) $, we have

\begin{displaymath}
\Vert D \vv \Vert _{L^{p}} ^{p}  =  \sum_{i \, j \,=\,1}^{3}  \int_{\mathbb{R}^{3}} | D_{j} v_{i}(x) |^{p} \, dx,
\end{displaymath}
and,  for general $ m \geq 1$
\begin{displaymath}
\Vert D^m \vv \Vert _{L^{p}} ^{p} =  \sum_{i\,=\,1}^{3} \sum_{j_1\,=\,1}^{3}  \sum_{j_2\,=\,1}^{3} \cdots  \sum_{j_m \,=\,1}^{3} \int_{\mathbb{R}^{3}} | D_{j_1} D_{j_2} \cdots D_{j_m} v_{i}(x) |^{p} \, dx,
\end{displaymath}
where $ D_j = \partial / \partial x_{j},   D_{j} D_{\ell} = \partial^{2} / \partial x_{j} \partial x_{\ell} $,  and so forth, while 
\begin{displaymath}
\Vert D \vv \Vert _{L^{\infty}} = \max \{ \esssup _{x \in \RR^3} |D_j v_i (x)|: 1 \leq j  \leq 3   \},
\end{displaymath}
and, more generally
\begin{displaymath}
\Vert D^m \vv \Vert _{L^{\infty}} = \max \{ \esssup _{x \in \RR^3} |D_{j_1} \cdots D_{j_m} v_i (x)|: 1 \leq j_1, \cdots, j_m  \leq 3   \},
\end{displaymath}
for $ m \geq 1 $,  if $p = \infty $.  In the text, we will only use $ p = 2 $ or $ p = \infty $.  We will also be using
the Sobolev space $ H^{m} $,  i.e., the space of all those functions in $L^{2}  $ whose derivatives of order $m$ are again in $ L^{2}$. If $ v = (v_1, v_2,  v_3) $, we write $\vv \in H^{m} $ when $ \vv_{i} \in H^{m} $ for all $1 \leq i \leq 3$. By $ L^{2}_{\sigma} $ we denote the space of vector functions $\vv = (v_1, v_2,  v_3) \in L^{2} $ with $ \nabla \cdot \vv = 0 $ in distributional sense,  that is,  $ \nabla \cdot \vv = 0 $ in $ {\mathcal{ D}}^{\prime} $ and similarly $ H^{m}_{\sigma}  = \{ \vv \in L^{2}_{\sigma}: \vv \in H^{m}\}$.  Here,   $ \nabla \cdot $ denotes the divergence operator; the curl and gradient are denoted  by $ \wedge$ and $ \nabla$ respectively.  Finally,  we will also occasionally mention  the homogeneous Sobolev space $ \dot{H}^{m},  m \geq 1 $, that is,  the space of tempered distributions with locally integrable Fourier transforms and whose distributional derivatives of order $m$  are all in $L^{2}$.  

\subsection{Review of results}
\label{review} 

In this section we present a brief review  of several basic results about (global) Leray solutions to  the micropolar systems \eqref{eqn:equation-u} - \eqref{eqn:equation-w} which will be needed later. The discussion below is adapted from Braz e Silva, J. Zingano and P. Zingano \cite{MR3907942},  Guterres, Melo, Nunes and Perusato \cite{MR3906315}, Guterres, Niche, Perusato and P. Zingano \cite{MR4546674}  and Hagstrom, Lorenz, J. Zingano and P. Zingano \cite{MR4021907}.  We consider general initial data $ \zz_0 \in L^{2}_{\sigma} \times L^{2}$.  We introduce the notation 

\begin{equation}
\label{eqn:gamma}
\gamma = \min \{ \mu,  \nu \},
\end{equation}
where $ \mu, \nu > 0 $ are given in \eqref{eqn:equation-u} and \eqref{eqn:equation-w}.  We also recall the general estimate 

\begin{equation}
\label{eqn:estimate-bzz}
\| D^{\ell} \vv(t) \|_{L^{\infty}} \, \| D^{m-\ell} \vv (t) \|_{L^{2}} \leq  \| \vv(t) \| ^{\frac{1}{2}} _{L^{2}} \, \Vert D \vv (t) \Vert ^{\frac{1}{2}} _{L^2} \, \| D^{m + 1} \vv(t) \|_{L^{2}}
\end{equation}
for every $ m \geq 1 $, $ 0 \leq \ell \leq m - 1 $ and $\vv \in H^{m+1}$,  see Braz e Silva, J. Zingano and P. Zingano \cite{MR3907942}. 

Uniqueness and exact regularity properties of solutions  are not known for arbitrary data, but the following result is available. 

\begin{Theorem}
\label{monotonicity}
Given $ \zz_0 \in L^{2}_{\sigma} \times L^{2}$ and Leray solutions $z$ to \eqref{eqn:equation-u} - \eqref{eqn:div-free},  there exists $t_{\ast \ast} \geq 0$ satisfying

\begin{displaymath}
t_{\ast\ast} \leq  0.005  \cdot \gamma^{- 5}  \, \| \zz_0 \| ^4 _{L^2}
\end{displaymath}
such that $\zz \in C^{\infty} \left( \RR^3 \times (t_{\ast \ast}, \infty) \right)$,  $\zz  \in C^{\infty}  \left( (t_{\ast \ast}, \infty), H^m \times H^m \right)$, for every $m \geq 0$ and $\Vert D \zz (t) \Vert _{L^2}$ is monotonically decreasing in $(t_{\ast \ast}, \infty)$. 
\end{Theorem}

{\bf Proof.} If $\zz \in C \ \left([t_0, T],  H^1 \right)$ for some $0 \leq t_0 < T$,  the solution is smooth in $\left( t_0, T \right)$ and differentiating \eqref{eqn:equation-u} - \eqref{eqn:equation-w} we obtain  the energy estimate 

\begin{align}
\label{eqn:first-estimate-monotonicity}
\Vert D \zz (t) \Vert _{L^2} ^2   & + 2 \gamma \int _{t_0} ^t \Vert D^2 \zz (\tau) \Vert ^2 _{L^2} \, d \tau \leq \Vert D \zz (t_0) \Vert ^2 _{L^2}  \notag \\ & + 2 \int _{t_0} ^t \Vert D^2 \zz (\tau) \Vert _{L^2} \Vert \zz (t) \Vert _{L^{\infty}} \Vert D \zz (\tau) \Vert _{L^2} \, d \tau 
\end{align}
for all $t \in \left( t_0, T \right)$, where $\gamma$ is as in \eqref{eqn:gamma}.  Note that

\begin{align}
\Vert \zz (\tau) \Vert _{L^{\infty}} \Vert D \zz (\tau) \Vert _{L^2}  & \leq K \, \Vert \zz (\tau) \Vert _{L^2} ^{\frac{1}{4}} \,  \Vert D \zz (\tau) \Vert _{L^2} \, \Vert D^2 \zz (\tau) \Vert )_{L^2} ^{\frac{3}{4}} \notag \\ & \leq K \, \Vert \zz (\tau) \Vert _{L^2} ^{\frac{1}{2}} \, \Vert D \zz (\tau) \Vert _{L^2} ^{\frac{1}{2}}  \, \Vert D^2 \zz (\tau) \Vert _{L^2},
\end{align}
where $ K =  \frac{\sqrt[8]{12}}{\sqrt{6 \pi}} $, see Theorem 2.2 in Schutz, Ziebell,  J. Zingano and P. Zingano \cite{MR3524859}. It then follows from \eqref{eqn:first-estimate-monotonicity} that $\Vert D \zz (t) \Vert _{L^2}$ is monotonically decreasing for $t \geq t_0$ if we have

\begin{equation}
\label{eqn:smallness}
K \, \Vert \zz (t) \Vert _{L^2} ^{\frac{1}{2}} \, \Vert D \zz (t) \Vert _{L^2} ^{\frac{1}{2}} \leq \gamma.
\end{equation}
In this case,  we have  $\zz \in C^{\infty} \left( \RR^3 \times (t_0, \infty) \right)$ and $\zz  \in C^{\infty}  \left( (t_0 \infty), H^m \times H^m \right)$ for all $m \geq 0$,  see Dong and Chen \cite{MR2572698},  Eringen \cite{MR0204005}, Galdi and Rionero \cite{MR0467030}. Taking $\hat{t} > \frac{1}{2} \, K^4 \, \gamma ^{ -5}  \,  \|  \zz_0 \| ^4 _{L^2}$,  from \eqref{eqn:energy-inequality} we obtain

\begin{displaymath}
2 \gamma \int _0 ^{\hat{t}} \Vert D \zz (\tau) \Vert _{L^2} ^2 \, d\tau \leq \Vert \zz_0 \| ^2 _{L^2},
\end{displaymath}
so there exists a set $E \subset (0,  \hat{t})$ with positive Lebesgue measure such that for all $t' \in E$

\begin{displaymath}
\Vert D \zz (t') \Vert _{L^2} ^2 \leq \frac{1}{2 \gamma} \,  \hat{t} ^{-1} \, \Vert  \zz_0 \| ^2 _{L^2}. 
\end{displaymath}
Therefore, for each $t' \in E$

\begin{displaymath}
\Vert \zz (t') \Vert _{L^2} ^2 \, \Vert D \zz (t') \Vert _{L^2} ^2  \leq \Vert \zz (t') \Vert _{L^2} ^2 \,  \frac{1}{2 \gamma} \, \hat{t} ^{-1} \, \Vert \zz_0 \| ^2 _{L^2}  \leq  \frac{1}{2 \gamma} \, \hat{t} ^{-1} \, \Vert \zz_0 \| ^4 _{L^2} \leq K ^{-4} \gamma ^4,
\end{displaymath}
where the last estimate follows from the choice of $\hat{t}$. The result follows from \eqref{eqn:smallness}.  $\Box$ \\

The monotonicity of $\Vert D \zz (t) \Vert _{L^2} ^2$ has important consequences,  as illustrated by the following results, more examples are given in Section \ref{proofs} (see also Guterres, Niche, Perusato and P. Zingano \cite{MR4546674}). 

\begin{Theorem}
\label{decay-dz}
For any Leray solution $\zz$ to \eqref{eqn:equation-u} - \eqref{eqn:div-free}, we have

\begin{equation}
\label{eqn:weak-decay-dz}
\Vert D \zz (t) \Vert _{L^2}  \leq \gamma ^{-\frac{1}{2}} \Vert \zz_0 \Vert _{L^2} \, t ^{-\frac{1}{2}},  \qquad t > 2 \, t_{\ast \ast},
\end{equation}
where $t_{\ast \ast}$ is the monotonicity time \eqref{eqn:monotonicity-time} of $\Vert D \zz (t) \Vert _{L^2}$.
\end{Theorem}

{\bf Proof.} From the energy inequality \eqref{eqn:energy-inequality} we have

\begin{displaymath}
\Vert \zz (t) \Vert _{L^2} ^2 + 2 \gamma \int_{t_0} ^t \Vert D \zz (\tau) \Vert _{L^2} ^2 \, d \tau \leq \Vert \zz _0 \Vert _{L^2} ^2,
\end{displaymath}
for every $t > t_0 > t_{\ast}$, where $t_{\ast}$ is the eventual regularity time of the solution. Then, for $t > 2 t_{\ast \ast}$ we have

\begin{displaymath}
\gamma \, t \, \Vert D \zz (t) \Vert _{L^2} ^2 \leq 2 \gamma \int_{\frac{t}{2}} ^t \Vert D \zz (\tau) \Vert _{L^2} ^2 \, d \tau \leq \Vert \zz_0 \Vert _{L^2} ^2. \qquad \Box
\end{displaymath}

\begin{Remark}
In a similar way we obtain that $\Vert D \zz (t) \Vert _{L^2} = o \left( t^{- \frac{1}{2}} \right)$, as for $t > t_{\ast \ast}$ we have

\begin{displaymath}
\frac{t}{2} \, \Vert D \zz (t) \Vert _{L^2} ^2 \leq \int_{\frac{t}{2}} ^{t} \Vert D \zz (\tau) \Vert _{L^2} ^2 \, d \tau \xrightarrow{t \to \infty} 0,
\end{displaymath} 
see Kreiss, Hagstrom, Lorenz and P. Zingano \cite{MR1994780}, Zhou \cite{MR2329641}. 
\end{Remark}

\begin{Theorem}
\label{monotonicity-dm}
For each $m \geq 0$ there exists a constant $K_m > 0$,  which depends on $m$ only,  such that for any Leray solution $\zz $ of \eqref{eqn:equation-u} - \eqref{eqn:div-free},  we have that $\Vert D ^m \zz (t) \Vert _{L^2}$ is monotonically decreasing in the interval $\left( t^{(m)} _{\ast \ast}, \infty \right)$ for some $ t^{(m)} _{\ast \ast} \geq t_{\ast}$  satisfying

\begin{equation}
\label{eqn:monotonicity-time-dm}
t_{\ast\ast} ^{(m)} \leq  K_m \cdot \gamma^{- 5} \,  \|  \zz_0 \| ^4 _{L^2}.
\end{equation}
\end{Theorem}

{\bf Proof.} The case $m \leq 1$ has already been considered.  Now,  for $m \geq 1$,  from    \eqref{eqn:equation-u} - \eqref{eqn:equation-w} we obtain 

\begin{align*}
\Vert D^m \zz (t) \Vert _{L^2} ^2 & + 2 \gamma \int_{t_0} ^t \Vert D^{m+1} \zz (\tau) \Vert _{L^2} ^2 \, d \tau \leq \Vert D^m \zz (t_0) \Vert _{L^2} ^2  \\ & + H'_m \sum _{\ell = 0} ^{[m/2]} \int _{t_0} ^t \Vert D^{m+1} \zz (\tau) \Vert_{L^2} \, \Vert D ^{\ell} \zz (\tau) \Vert _{L^{\infty}} \, \Vert D^{m - \ell} \zz (\tau) \Vert_{L^2}   \, d \tau 
\end{align*}
for all $t> t_0 > t_{\ast}$,  where $[r]$ denotes the integer part of $r \in \RR$ and $H'_m >0$ is a constant that depends on $m$ only.   From this and \eqref{eqn:estimate-bzz}we obtain

\begin{align}
\label{eqn:condition-monotonicity-dm}
\Vert D^m \zz (t) \Vert _{L^2} ^2 & + 2 \gamma \int_{t_0} ^t \Vert D^{m+1} \zz (\tau) \Vert _{L^2} ^2 \, d \tau \leq \Vert D^m \zz (t_0) \Vert _{L^2} ^2 \notag  \\ & + H_m  \int _{t_0} ^t g (\tau) \, \Vert D^{m+1} \zz (\tau) \Vert_{L^2}  \, d \tau 
\end{align}
where $g(\tau) = \Vert \zz (\tau) \Vert _{L^2} ^{\frac{1}{2}} \, \Vert D \zz (\tau) \Vert _{L^2} ^{\frac{1}{2}}$.   From \eqref{eqn:condition-monotonicity-dm}, we see that monotonicity of $\Vert D ^m \zz (t) \Vert _{L^2}$ is achieved if $H_m \, g(t) \leq 2 \gamma$, when $t > t^{(m)} _{\ast \ast}$.  From \eqref{eqn:monotonicity-time}, \eqref{eqn:weak-decay-dz} and \eqref{eqn:monotonicity-time-dm} , this condition is guaranteed if $K_m = \frac{H_m ^4}{16}$. $\Box$ \\

\begin{Remark}
By slightly redefining the values of $K_m$ and $t_{\ast \ast} ^{(m)}$ we can obtain more information about  $\Vert D^m \zz (t) \Vert _{L^2}$.  By setting $K_m = H_m ^4$, for $m \geq 1$, we obtain that $H_m \, g_n(t) \leq \gamma$, for $t > t_{\ast \ast} ^{(m)}$, where $t > t_{\ast \ast} ^{(m)}$ has been redefined as

\begin{displaymath}
t > t_{\ast \ast} ^{(m)} = H_m ^4 \, \gamma ^{-5} \, \Vert \zz_0 \Vert _{L^2} ^4.
\end{displaymath}
Then, \eqref{eqn:condition-monotonicity-dm} leads to

\begin{equation}
\label{eqn:-estimate-dmz-basic}
\Vert D^m \zz (t) \Vert _{L^2} ^2  + \gamma \int_{t_0} ^t \Vert D^{m+1} \zz (\tau) \Vert _{L^2} ^2 \, d \tau \leq \Vert D^m \zz (t_0) \Vert _{L^2} ^2.
\end{equation}
\end{Remark}

\begin{Theorem}
\label{decay-dmz}
For any Leray solution $\zz$ to \eqref{eqn:equation-u} - \eqref{eqn:div-free}, we have

\begin{equation}
\label{eqn:weak-decay-dmz}
\Vert D ^m \zz (t) \Vert _{L^2}  \leq C  \, t ^{-\frac{m}{2}},  \qquad t > 2 \, t_{\ast \ast} ^{(m)},
\end{equation}
where $t_{\ast \ast} ^{(m)}$ is the monotonicity time \eqref{eqn:monotonicity-time-dm} of $\Vert D ^m \zz (t) \Vert _{L^2} ^2$ and the constant $C = C \left(m, C_0,  \gamma \right)$.
\end{Theorem}

{\bf Proof.}  The case $m = 1$ is just Theorem \ref{decay-dz}.  In order to prove the estimate for $m = 2$,  take $m = 1$ in \eqref{eqn:-estimate-dmz-basic}.  Then, for $t > 2 t_{\ast \ast} ^{(2)}$ we have

\begin{displaymath}
\gamma \,  \frac{t}{2} \, \Vert D ^2 \zz (t) \Vert _{L^2} ^2 \leq \gamma \int _{\frac{t}{2}} ^t \Vert D ^2 \zz (\tau) \Vert _{L^2} ^2 \, d \tau \leq \Vert D \zz (t/2) \Vert _{L^2} ^2, 
\end{displaymath}
where the first inequality is due to the monotonicity of $\Vert D ^2 \zz (\tau) \Vert _{L^2} ^2$,  as $\frac{t}{2} > t_{\ast \ast} ^{ (2)}$.  For $t > 4 t _{\ast \ast} ^{(2)}$, we have $\frac{t}{2} > 2 t_{\ast \ast} ^{(1)}$,  because  $t _{\ast \ast} ^{(2)} > t _{\ast \ast} ^{(1)}$.  Then,  from the previous estimate we deduce

\begin{displaymath}
\Vert D \zz (t/2) \Vert _{L^2} ^2 \leq K_1 ^2 \, \gamma ^{-1} \Vert \zz_0 \Vert _{L^2} ^2 \, \left( \frac{t}{2} \right) ^{-1},
\end{displaymath}
which proves our estimate.  An induction procedure proves the result for any $m$.  $\Box$ \\

\begin{Theorem}
\label{decay-dmw}
For any Leray solution $\zz = (\uu,  \ww)$ to \eqref{eqn:equation-u} - \eqref{eqn:div-free}, we have

\begin{equation}
\label{eqn:weak-decay-dw}
\Vert \ww (t) \Vert _{L^2} \leq C \, t ^{-\frac{1}{2}},  \qquad t > \tau_0,
\end{equation}
for some $\tau_0  = \tau_0  (\chi, \gamma, t_{\ast \ast}) > 0$ and $C = C (\gamma, \Vert \zz_0 \Vert _{L^2})$.  Moreover,  for every $m \geq 1$ we have

\begin{equation}
\label{eqn:weak-decay-dmw}
\Vert D^m \ww (t) \Vert _{L^2}  \leq C \, t ^{-\frac{m+1}{2}},  \qquad t > \tau_m,
\end{equation}
for some constant $C = C \left(m, C_0,  \gamma \right)   > 0$  and $\tau_m = \tau_m \left(m , \gamma, \chi, t_{\ast \ast} ^{(m+2)},  \Vert \zz_0 \Vert _{L^2} \right) > 0$. 
\end{Theorem}

{\bf Proof.} From \eqref{eqn:equation-w} we obtain, for every $t > t_{\ast}$

\begin{align*}
\frac{d}{dt} \Vert \ww (t) \Vert _{L^2} ^2 & + 2 \nu \, \Vert D \ww (t) \Vert _{L^2} ^2 + 8 \chi \, \Vert \ww(t) \Vert _{L^2} ^2  \leq 4 \chi \, \Vert \ww (t) \Vert _{L^2} \, \Vert D \uu (t) \Vert _{L^2}  \\ \leq \, & 4 \chi \, \Vert \ww (t) \Vert _{L^2} ^2 + \chi \, \Vert D \uu (t) \Vert _{L^2} ^2.
\end{align*}
Then

\begin{displaymath}
\frac{d}{dt} \Vert \ww (t) \Vert _{L^2} ^2 + 4 \chi \, \Vert \ww(t) \Vert _{L^2} ^2 \leq \chi \, \Vert D \uu (t) \Vert _{L^2} ^2.
\end{displaymath}
Using \eqref{eqn:weak-decay-dz}, we obtain that for $t > 4 t_{\ast \ast}$

\begin{align*}
\Vert \ww (t) \Vert _{L^2} ^2 & \leq  e ^{-2 \chi t} \, \Vert \zz_0 \Vert _{L^2} ^2 + \chi \int _{\frac{t}{2}} ^t e^{-4 \chi (t - \tau)} \, \Vert D \uu (\tau) \Vert _{L^2} ^2 \, d \tau \\ & \leq \Vert \zz_0 \Vert _{L^2} ^2 \left( e^{-2 \chi t} + \frac{1}{2} \,  \gamma^{-1} \, t^{-1}  \right),
\end{align*}
which,  for large enough $t$,  leads to \eqref{eqn:weak-decay-dw}.

Now,  differentiating \eqref{eqn:equation-w} $m$ times, we obtain for $t > t_{\ast}$

\begin{align*}
\frac{d}{dt} \Vert D^m \ww (t) \Vert _{L^2} ^2 & + 2 \nu \, \Vert D ^{ m+1} \ww (t) \Vert _{L^2} ^2 + 8 \chi \, \Vert D^m \ww(t) \Vert _{L^2} ^2  \\ & \leq 4 \chi \, \Vert D^m \ww (t) \Vert _{L^2} \, \Vert D ^{m+1}  \uu (t) \Vert _{L^2}  + K_m ' \Vert D^m \ww (t) \Vert _{L^2} \, G_m (t)
\end{align*}
where $ G_{m}(t) = \sum_{\ell = 0}^{m} \Vert D^{\ell} \uu ( t) \Vert _{L^{\infty}} \Vert D^{m+1 -\ell} \ww (t) \Vert _{L^2}$,  for some constant $K_m ' > 0$, which depends on $m$ only.  This leads to

\begin{displaymath}
\frac{d}{dt} \Vert D^m \ww (t) \Vert _{L^2} ^2  + 2 \chi \, \Vert D ^{ m} \ww (t) \Vert _{L^2} ^2 \leq \chi \, \Vert D ^{m+1}  \uu (t) \Vert _{L^2} ^2  + K_m '' \, \chi ^{-1} \, \Vert D^m \ww (t) \Vert _{L^2} \, G_m ^2 (t)
\end{displaymath}
where $K_m '' = \left( \frac{K_m '}{\sqrt{8}} \right) ^2$. Then, for $t > 2 t_{\ast}$ we have

\begin{align*}
\Vert D^m \ww (t) \Vert _{L^2} ^2 & \leq  e ^{- \chi t} \, \Vert D^m \ww (t/2) \Vert _{L^2} ^2 + \chi \int _{\frac{t}{2}} ^t e^{-2 \chi (t - \tau)} \, \Vert D ^{m+1} \uu (\tau) \Vert _{L^2} ^2 \, d \tau \\ & + K_m ''\, \chi ^{-1} \int_{\frac{t}{2}} ^t  e^{-2 \chi (t - \tau)} \, G_m ^2 (\tau) \, d \tau.
\end{align*}
Taking $\frac{t}{2} > 2^{(m+2)} \, t_{\ast\ast}^{(m+2)}$ and using \eqref{eqn:weak-decay-dmz} to  to estimate the three terms on the righthand side of the expression above, we obtain \eqref{eqn:weak-decay-dmw}.  $\Box$ \\ 

\begin{Remark}
 \label{bootstrap-dmz-dmw}
Proceeding as in the derivation of \eqref{eqn:weak-decay-dz} and \eqref{eqn:weak-decay-dmz}, we obtain from \eqref{eqn:weak-decay-dz} and \eqref{eqn:-estimate-dmz-basic} that

\begin{displaymath}
\Vert D \zz (t) \Vert _{L^2} \leq \gamma ^{- \frac{1}{2}} \, \Vert \zz (t/2) \Vert _{L^2} \, t ^{-\frac{1}{2}}, \qquad t > 2 t _{\ast \ast},
\end{displaymath}
where $t_{\ast \ast}$ is the monotonicity time \eqref{eqn:monotonicity-time}.  Also, for every $m > 1$

\begin{equation}
\label{eqn:estimate-dmz-bootstrap}
\Vert D^m \zz (t) \Vert _{L^2} \leq \widetilde{K} _m \, \gamma ^{- \frac{m}{2}} \, \Vert \zz (t/2) \Vert _{L^2} \, t ^{- \frac{m}{2}}, \qquad t > 2^m t_{\ast \ast} ^{(m)},
\end{equation}
for some constant $\widetilde{K}_m$ which depends on $m$ only and $t_{\ast \ast} ^{(m)}$ as in \eqref{eqn:monotonicity-time-dm}. Through these estimates we can obtain, as in Theorem \ref{decay-dmw}, that for every $m \geq 0$

\begin{equation}
\label{eqn:estimate-dmw-bootstrap}
\Vert D ^m \ww (t) \Vert _{L^2} \leq \widetilde{C} _m \, \gamma ^{- \frac{m+1}{2}} \, \Vert \zz (t/2) \Vert _{L^2} \, t^{- \frac{m+1}{2}}, \qquad t > \widetilde{\tau}_m 
\end{equation}
where $\widetilde{C}_m$ depends on $m$ only  and $\widetilde{\tau}_m$ behaves like $\tau_m$ in \eqref{eqn:weak-decay-dmw}. 
\end{Remark}

\section{Proof of Theorems \ref{decay-du-dw-upper},  \ref{decay-error-upper} and \ref{lower-bound-du}} 
\label{proofs}

In this section we derive Theorems \ref{decay-du-dw-upper}, \ref{decay-error-upper} and \ref{lower-bound-du}, dividing their proofs into several individual results of interest on their own.  Throughout the section $\zz = \left(\uu,   \ww  \right)$ denotes a Leray solution to \eqref{eqn:equation-u} - \eqref{eqn:equation-w} with initial data $\zz_0 = \left(\uu_0, \ww_0 \right) \in L^2 _{\sigma} \times L^2$,  for which \eqref{eqn:upper-bound-u} holds, i.e.  for some $\alpha > 0$ we have

\begin{displaymath}
\Vert \uu (t) \Vert _{L^2} \leq C t ^{- \alpha}, \qquad \forall t \geq t_0.
\end{displaymath}
We also recall that $\gamma = \min \{\mu, \nu \}$. 

\subsection{Preliminary results} We first prove some preliminary decay estimates which will be used in the proof of the main results stated in Section \ref{main-results}.  Our first result is a version of Theorem \ref{decay-du-dw-upper},  but with the constant depending on $\gamma = \min \{\mu, \nu \}$ instead of just $\mu$.

\begin{Theorem}
\label{decay-du-dw-upper-gamma}
If \eqref{eqn:upper-bound-u} is valid for some $\alpha \geq 0$,  then we have that for $m \geq 1$

\begin{equation}
%\label{eqn:upper-bound-du}
\Vert D^m \uu(t) \Vert _{L^2}  \leq C \, \gamma ^{- \frac{m}{2}} \, t ^{- \left( \alpha + \frac{m}{2} \right)},  \qquad \forall \, t \geq T_m,
\end{equation}
where $C = C (C_0, m, \alpha)$ and $T_m >  t_{\ast \ast}$,  see Remark \ref{dependence-times-constants} for the dependence of $T_m$ on the constants appearing in the problem.  Moreover, we have that for every $m \geq 0$

\begin{equation}
\label{eqn:upper-bound-dw-three-two}
\Vert D^m \ww(t) \Vert _{L^2} \leq  C \, \gamma ^{- \frac{m + 1}{2}} \,  t ^{- \left( \alpha + \frac{m + 1}{2} \right)},  \qquad \forall \, t \geq T'_m,
\end{equation}
where $C = C (C_0, m+1, \alpha)$ and  $T'_m >  t_{\ast \ast}$,  see Remark \ref{dependence-times-constants}. 
\end{Theorem}

{\bf Proof:} From \eqref{eqn:estimate-dmz-bootstrap} and \eqref{eqn:estimate-dmw-bootstrap} in Remark \ref{bootstrap-dmz-dmw} and \eqref{eqn:upper-bound-u}, we only need to prove that

\begin{displaymath}
\Vert \ww (t) \Vert _{L^2} \leq C t^{- \alpha}, \qquad t \geq T_0 '
\end{displaymath}
for some large enough $T_0 ' > 0$ and $C = C \left(\alpha,  \Vert \zz _0 \Vert _{L^2} \right)$.  From \eqref{eqn:equation-w} we obtain

\begin{displaymath}
\frac{d}{dt} \Vert \ww (t) \Vert _{L^2} ^2 + 2 \nu \Vert D \ww (t) \Vert _{L^2} ^2 + 8 \chi \Vert \ww (t) \Vert _{L^2} ^2 \leq 4 \chi \Vert D \ww (t) \Vert _{L^2} \Vert u(t) \Vert _{L^2}, 
\end{displaymath}
for $t > t_{\ast}$, where $t_{\ast}$ is the solution's eventual regularity time.  This implies that

\begin{displaymath}
\frac{d}{dt} \Vert \ww (t) \Vert _{L^2} ^2  + 8 \chi \Vert \ww (t) \Vert _{L^2} ^2 \leq 2 \nu ^{ -1} \chi ^2 \Vert u(t) \Vert ^2 _{L^2}, 
\end{displaymath}
 which leads,  through \eqref{eqn:upper-bound-u} and for $t_1 > \frac{t}{2} > \max \{t_{\ast}, T_0 \}$,  to 

\begin{align*}
\Vert \ww (t) \Vert _{L^2} ^2 & \leq C e^{- 4 \chi t} + 2 \nu ^{ -1} \chi ^2 \int _{\frac{t}{2}} ^t e ^{-8 \chi  (t - \tau)} \Vert \uu (\tau) \Vert _{L^2} ^2 \, d \tau \\ & \leq C e^{- 4 \chi t} + C_1 t^{-2 \alpha},  \qquad t > t_1
\end{align*}
where $C = C (\Vert \zz_0 \Vert_{L^2} ^2)$ and $C_1 = C_1 \left(\alpha, \nu, \chi \right) $.  Increasing $t_1$ if necessary, 
the result for $\Vert \ww (t) \Vert _{L^2}$ is obtained.  $\Box$

We now estimate the ``error'' \eqref{eqn:error},  i.e.

\begin{displaymath}
\error (x,t) = \ww(x,t) -  \frac{1}{2} \nabla \wedge \uu (x,t),
\end{displaymath}
for $t > t_{\ast}$, where $t_{\ast}$ is the solution's regularity time.  From \eqref{eqn:equation-u} and \eqref{eqn:equation-w} we obtain

\begin{equation}
\label{eqn:equation-error}
\partial_t \error + \left(\uu \cdot \nabla \right) \error + 4 \chi \error = \mathbb{L} \error + (\nu - \mu) \Delta \ww + \frac{1}{2} \sum _{j=1} ^3 \left(\nabla u_j \right) \wedge  \left(\nabla  D_j u \right),
\end{equation}
where $\mathbb{L}$ denotes the elliptic operator

\begin{equation}
\label{eqn:operator-l}
\mathbb{L} \vv = \mu \Delta  \vv + \kappa \nabla \left( \nabla \cdot \vv \right) - \chi \nabla \wedge \left( \nabla \wedge \vv \right).
\end{equation}
The following preliminary estimate is key for proving \eqref{eqn:upper-bound-error} in Theorem \ref{decay-error-upper}.

\begin{Theorem}
\label{estimate-error-first}
Suppose \eqref{eqn:upper-bound-u} holds.  Then, for any $m \geq 0$

\begin{equation}
\label{eqn:first-estimate-dmerror}
\Vert D^m \error (t) \Vert _{L^2} \leq C\, \gamma^{- \frac{m+3}{2}} \,  t ^{- \left( \alpha + \frac{m+3}{2} \right)}, \qquad t > T_m '',
\end{equation}
where $C = C(m, \alpha, C_0, \mu, \nu, \chi)$ and $T_m ''$ depends on the constants in the problem as explained in Remark \ref{dependence-times-constants}.
\end{Theorem}

{\bf Proof:} Estimate \eqref{eqn:first-estimate-dmerror} will follow from  \eqref{eqn:upper-bound-u} and 

\begin{equation}
\Vert D^m \error (t) \Vert _{L^2} \leq C \Vert \zz (t/2) \Vert _{L^2} t ^{-\frac{m+3}{2}}, \qquad t > \sigma _m,
\end{equation}
for some $\sigma_m$ that depends on the constants in the problem as stated in Remark \ref{dependence-times-constants}.  From \eqref{eqn:equation-error} and \eqref{eqn:operator-l} we obtain

\begin{align*}
\frac{d}{dt} \Vert D^m \error (t) \Vert _{L^2} ^2 & + 2 \mu \Vert D^{m+1} \error (t) \Vert _{L^2} ^2 + 8 \chi \Vert D^m \error (t)  \Vert _{L^2} ^2 \\ & \leq |\mu - \nu| \Vert D^m \error (t) \Vert _{L^2}  \Vert D^{m+2} \ww (t) \Vert _{L^2}  + K_m \Vert D^m \error (t)  \Vert _{L^2} \, H_m (t),   
\end{align*}
for $t > t_{\ast}$,  where

\begin{displaymath}
H_m (t) = \sum _{\ell = 0} ^m \Vert D^{\ell} \uu (t) \Vert _{L^{\infty}} \Vert D^{m+1-\ell} \error (t)  \Vert _{L^2} + \sum _{\ell = 0} ^m \Vert D^{\ell + 1} \uu (t) \Vert _{L^{\infty}} \Vert D^{m+1-\ell} \uu (t)  \Vert _{L^2},
\end{displaymath}
for some constant $K_m  > 0$ that depends on $m$ only.  From this we obtain

\begin{equation}
\label{eqn:intermediate-estimate-dmerror}
\frac{d}{dt} \Vert D^m \error (t) \Vert _{L^2} ^2  + 4 \chi \Vert D^m \error (t)  \Vert _{L^2} ^2 \leq C \Vert D^{m+2} \ww (t) \Vert _{L^2} ^2 + C H_m ^2 (t).
\end{equation}
The expression for $H_m$ depends on $D^{m+1-\ell} \error (t)$,  which can in turn be expressed through $D^{m+1-\ell} \ww (t)$ and $D^{m+2-\ell} \uu (t)$, because of the definition for $\error$ in \eqref{eqn:error}.  And then,  the norms

\begin{align*}
\Vert D^{m+2} \ww (t)  \Vert _{L^2} &,  \Vert D^{\ell} \uu (t) \Vert _{L^{\infty}},  \Vert D^{m+1-\ell} \ww (t) \Vert _{L^2},  \Vert D^{m+2-\ell} \uu (t)\Vert _{L^2} \\ & \Vert D^{\ell + 1} \uu (t) \Vert _{L^{\infty}},   \Vert D^{m+1-\ell} \uu (t)  \Vert _{L^2},
\end{align*}
on the right hand side of \eqref{eqn:intermediate-estimate-dmerror} can be estimated through \eqref{eqn:estimate-dmz-bootstrap} and \eqref{eqn:estimate-dmw-bootstrap} in Remark \ref{bootstrap-dmz-dmw},  interpolation and \eqref{eqn:upper-bound-dw-three-two}.    After this,  integrating \eqref{eqn:intermediate-estimate-dmerror} and dropping the exponential term,  we obtain the result.  $\Box$

\subsection{Proofs of Main Results} In order to prove Theorems \ref{decay-du-dw-upper},  \ref{decay-error-upper} and \ref{lower-bound-du} using the preliminary estimates from the previous Section, we rewrite \eqref{eqn:equation-u} and \eqref{eqn:equation-w} through \eqref{eqn:error}, as

\begin{align}
\label{eqn:equation-u-rewritten}
\partial_t \uu + \uu \cdot \nabla \uu + \nabla p & =  \mu  \Delta \uu + 2 \chi \nabla \wedge \error,  \\
\label{eqn:equation-w-rewritten}
\partial_t \ww + \uu \cdot \nabla \ww & =  \nu \Delta \ww + \kappa \nabla (\nabla \cdot \ww) - 4 \chi \error,  \\
\div \, \uu & =  0.
\end{align}
In this new setting,  we first obtain a refined monotonicity estimate.

\begin{Theorem}
\label{monotonicity-zeta}
Suppose \eqref{eqn:upper-bound-u} holds. Then, for any $m \geq 1$ the function $Z_m(t)$ given by

\begin{equation}
\label{eqn:zeta-m}
Z_m (t) = \Vert D^m \uu (t) \Vert _{L^2} ^2 + 4 \chi^2 \mu ^{-1} \int _{t} ^{\infty} \Vert D^m \error (\tau) \Vert _{L^2} ^2 \, d \tau,
\end{equation}
is monotonically decreasing in an interval $\left(\zeta_m, \infty \right)$, where $\zeta_m$ depends on the constants described in Remark \ref{dependence-times-constants}.
\end{Theorem}

{\bf Proof of Theorem \ref{monotonicity-zeta}:} From \eqref{eqn:equation-u-rewritten}, through \eqref{eqn:estimate-bzz} we obtain 

\begin{displaymath}
\frac{d}{dt} \Vert D^m \uu (t) \Vert _{L^2} ^2 + \mu \Vert D^{m+1} \uu (t) \Vert _{L^2} ^2 \leq 4 \chi \Vert D^{m+1} \uu (t) \Vert _{L^2} \Vert D^m \error (t) \Vert _{L^2}, 
\end{displaymath}
for $t > \zeta_m$,  for a large enough $\zeta_m$ which depends on the constants described in Remark \ref{dependence-times-constants}.  From this we deduce that

\begin{displaymath}
\frac{d}{dt} \Vert D^m \uu (t) \Vert _{L^2} ^2 \leq 4 \chi^2 \mu ^{ -1}  \Vert D^m \error (t) \Vert _{L^2} ^2, 
\end{displaymath}
from which the result follows.  $\Box$

\begin{Remark}
By a similar and simpler argument, it follows directly from \eqref{eqn:equation-u-rewritten} that

\begin{displaymath}
Z_0 (t) = \Vert \uu (t) \Vert _{L^2} ^2 + 2 \chi^2 \mu ^{-1} \int _{t} ^{\infty} \Vert \error (\tau) \Vert _{L^2} ^2 \, d \tau,
\end{displaymath}
decreases monotonically in $\left(t_{\ast} \infty \right)$, where $t_{\ast}$ is the regularity time.  We will not use this estimate in this work,  but we emphasize its importance,  as well as that of \eqref{eqn:zeta-m}, see the general discussion in Guterres, Niche, Perusato and P. Zingano \cite{MR4546674}. 
\end{Remark}

Now we can prove our first main result. 

{\bf Proof of Theorem \ref{decay-du-dw-upper}:} For $m = 1$,  from \eqref{eqn:equation-u-rewritten} we obtain, for $t > 2 t_{\ast}$,  that 

\begin{displaymath}
\Vert \uu(t) \Vert _{L^2} ^2 + 2 \mu \int_{\frac{t}{2}} ^t \Vert D \uu (\tau) \Vert _{L^2} ^2 \, d \tau \leq \Vert \uu (t/2) \Vert _{L^2} ^2 + 4 \chi \int_{\frac{t}{2}} ^t \Vert D \uu (\tau) \Vert _{L^2} \Vert \error (\tau) \Vert _{L^2} \, d \tau. 
\end{displaymath}
We then have

\begin{displaymath}
\Vert \uu(t) \Vert _{L^2} ^2 +  \mu \int_{\frac{t}{2}} ^t \Vert D \uu (\tau) \Vert _{L^2} ^2 \, d \tau \leq \Vert \uu (t/2) \Vert _{L^2} ^2 + 4 \chi^2 \mu ^{-1}  \int_{\frac{t}{2}} ^t  \Vert \error (\tau) \Vert _{L^2} ^2 \, d \tau. 
\end{displaymath}
Using \eqref{eqn:zeta-m} from Theorem \ref{monotonicity-zeta} we deduce that

\begin{displaymath}
\mu \frac{t}{2} \Vert D \uu (t) \Vert _{L^2} ^2 \leq \mu \int _{\frac{t}{2}} ^t Z_1 (\tau) \, d \tau \leq \Vert \uu \left( t/2 \right) \Vert  _{L^2} ^2  + 4 \chi ^2 \mu ^{-1} \int _{\frac{t}{2}} ^t \left( \Vert \error (\tau) \Vert _{L^2} ^2 + E_1 (\tau)  \right) \, d \tau,
\end{displaymath}
for $t > 2 \zeta_1$ and

\begin{displaymath}
E_1 (\tau) =  \mu \int _{\tau} ^{\infty} \Vert D \error (s) \Vert _{L^2} ^2 \, ds. 
\end{displaymath}
Then,  \eqref{eqn:upper-bound-u} and \eqref{eqn:first-estimate-dmerror} in Theorem \ref{estimate-error-first} lead us to the result for $D \uu (t)$.   The general case $m > 1$ can be obtained by induction using a similar argument.  To wit, note that after taking  derivatives in \eqref{eqn:equation-u-rewritten} and using \eqref{eqn:estimate-bzz} we obtain

\begin{displaymath}
\mu \int_{\frac{t}{2}} ^t \Vert D ^m \uu (\tau) \Vert _{L^2} ^2 \, d \tau \leq \Vert D ^{m-1} \uu (t/2) \Vert _{L^2} ^2 + 8 \chi ^2 \mu ^{-1}  \int_{\frac{t}{2}} ^t  \Vert D ^{m-1} \error (\tau) \Vert _{L^2} ^2 \, d \tau,
\end{displaymath}
for $t > t_m$,  where $t_m$ depends on the constants in Remark \ref{dependence-times-constants}.  From  \eqref{eqn:zeta-m} from Theorem \ref{monotonicity-zeta} and taking a larger $t_m$,  if necessary,  we obtain

\begin{displaymath}
\mu \frac{t}{2} \Vert D ^m \uu (t) \Vert _{L^2} ^2  \leq \Vert D ^{m-1} u \left( t/2 \right) |\Vert _{L^2} ^2  + 8 \chi ^2 \mu ^{-1} \int _{\frac{t}{2}} ^t \left( \Vert D ^{m-1} \error (\tau) \Vert _{L^2} ^2 + E_m (\tau)  \right) \, d \tau,
\end{displaymath}
where 

\begin{displaymath}
E_m (\tau) = \mu \int _{\tau} ^{\infty} \Vert D ^m \error (s) \Vert _{L^2} ^2 \, ds. 
\end{displaymath}
Using \eqref{eqn:first-estimate-dmerror} in Theorem \ref{estimate-error-first} and the induction hypothesis for $\Vert D ^{m-1} u \left( t/2 \right) \Vert _{L^2} ^2$,  we obtain the estimate \eqref{eqn:upper-bound-du}. 

From \eqref{eqn:error} we obtain that

\begin{displaymath}
\Vert D^m \ww (t) \Vert _{L^2} \leq \frac{1}{2} \Vert D^{m+1} \uu (t) \Vert _{L^2} + \Vert D^m \error (t) \Vert _{L^2}, 
\end{displaymath}
for $t > t_ {\ast}$.  Then,  \eqref{eqn:upper-bound-dw} is an immediate consequence of \eqref{eqn:upper-bound-u} and \eqref{eqn:zeta-m}. $\Box$

We now prove the full decay estimate for $\error$ as in \eqref{eqn:upper-bound-error},  thus improving Theorem \ref{estimate-error-first}.

{\bf Proof of Theorem \ref{decay-error-upper}:} Let $t_{\ast}$ be the eventual regularity time.  From  \eqref{eqn:equation-error} we deduce  that for $t > t_{\ast}$

\begin{align*}
\frac{d}{dt} \Vert D^m \error (t) \Vert _{L^2} ^2 & + 2 \mu \Vert D^{m+1} \error (t) \Vert _{L^2} ^2 + 8 \chi \Vert D^m \error (t)  \Vert _{L^2} ^2 \\  & \leq \Vert D^m \error (t) \Vert _{L^2} \left(2 |\mu - \nu| \Vert D^{m+2} \ww (t) \Vert _{L^2} + A_m F_m (t) + B_m G_m (t) \right),
\end{align*}
where $A_m, M_m$ are constants that depend only on $m$,  and

\begin{align*}
F_m(t) & = \sum _{\ell = 0} ^m \Vert D^{\ell + 1} \uu (t)  \Vert _{L^{\infty}} \Vert D^{m - \ell + 1} \uu (t) \Vert _{L^2},  \\ G_m(t) & = \sum _{\ell = 0} ^m \Vert D^{\ell} \uu (t)  \Vert _{L^{\infty}} \Vert D^{m - \ell + 1} \error (t) \Vert _{L^2}.
\end{align*}
From this we obtain 

\begin{displaymath}
\frac{d}{dt} \Vert D^m \error (t) \Vert _{L^2} ^2  + 4 \chi \Vert D^m \error (t) \Vert _{L^2} ^2 \leq C |\mu - \nu|^2 \Vert D^{m+2} \ww (t) \Vert _{L^2} ^2 + C F_m ^2 (t) + C G_m ^2 (t).
\end{displaymath}
Then,  using \eqref{eqn:first-estimate-dmerror}, \eqref{eqn:upper-bound-du} and \eqref{eqn:upper-bound-dw} in the expressions for $F_m(t), G_m (t)$ and $D^{m+2} \ww(t)$ and integrating in time the previous line,  we obtain  estimate \eqref{eqn:upper-bound-error}.  $\Box$

\begin{Remark}
From \eqref{eqn:error},  we have that $\nabla \cdot \ww = \nabla \cdot \error$.  Then,  we deduce that $\nabla \cdot \ww$ has faster decay than $\nabla \wedge \ww$,  as from \eqref{eqn:upper-bound-error} in Theorem \ref{decay-error-upper} we have

\begin{displaymath}
\Vert D^m \left( \nabla \cdot \ww (t) \right) \Vert _{L^2} \leq C_1 \, |\mu - \nu| \, t ^{- \left( \alpha + \frac{m + 4}{2}\right)} + C_2 \, t^{- \left( 2 \alpha + \frac{m+4}{2} + \frac{1}{4} \right)}, 
\end{displaymath}
while \eqref{eqn:upper-bound-dw} in Theorem \ref{decay-du-dw-upper} leads to
\begin{displaymath}
\Vert D^m \left(\nabla  \wedge \ww(t) \right) \Vert _{L^2} \leq  C \,   t ^{- \left( \alpha + \frac{m + 2}{2} \right)}.
\end{displaymath}
\end{Remark}

We prove now the lower bounds for $D^m \uu$.

{\bf Proof of Theorem \ref{lower-bound-du}:} multiplying by $\uu$ and integrating in \eqref{eqn:equation-u-rewritten}, we obtain

\begin{displaymath}
\Vert \uu (T) \Vert _{L^2} ^2 + 2 \mu \int_t ^T \Vert D \uu (\tau) \Vert _{L^2} ^2 \, d \tau= \Vert \uu (t) \Vert _{L^2} ^2 + 4 \chi \int _t ^T \langle \nabla \wedge \uu,  \error \rangle _{L^2} \, d \tau,
\end{displaymath}
for all $T > t > \max \{t_0, T_0, t_{\ast} \}$,  for $t_{\ast}$ the eventual regularity time  and $T_0, t_0$ as in \eqref{eqn:upper-bound-u} and \eqref{eqn:lower-bound-u} respectively.  Using these estimates and setting $\lambda_1 = \frac{C_0}{\widetilde{C}_0} \geq 1$, we obtain

\begin{align*}
4 \mu \int_t ^T \Vert D \uu (\tau) \Vert _{L^2} ^2 \, d \tau &  \geq \Vert \uu(t) \Vert _{L^2} ^2 - \Vert \uu(T) \Vert _{L^2} ^2 - 2 \chi ^2 \mu ^{-1} \int _t ^T \Vert \error (\tau) \Vert _{L^2} ^2 \, d \tau \\ & \geq C \left(t^{- 2 \alpha} - \lambda_1 ^2 \, T^{-2 \alpha} \right) - 2 \chi ^2 \mu ^{-1} \int _t ^T \Vert \error (\tau) \Vert _{L^2} ^2 \, d \tau.
\end{align*}
Taking $T = \left(2 \lambda_1 \right) ^{\frac{1}{\alpha}}  t$ and using \eqref{eqn:zeta-m} in Theorem \ref{monotonicity-zeta} we obtain, after choosing $t$ large enough

\begin{displaymath}
4 \mu  \left(2 \lambda_1 \right) ^{\frac{1}{\alpha}}  t \Vert D \uu (t) \Vert _{L^2} ^2 \geq C \, t^{-2 \alpha}  - 2 \chi ^2 \mu ^{-1} \int _t ^{\infty} \Vert \error (\tau) \Vert _{L^2} ^2 \, d \tau - C \int_t ^{\infty} \Vert D \error (\tau) \Vert _{L^2} ^2 \, d \tau,
\end{displaymath}
from which we deduce our estimate for $m = 1$.  We now use induction to prove the result.  From  \eqref{eqn:equation-u-rewritten} and \eqref{eqn:estimate-bzz} 

\begin{align*}
4 \mu  \int_t ^T \Vert D^m \uu (\tau) \Vert _{L^2} ^2 \, d \tau &  \geq \Vert D ^{m-1} \uu(t) \Vert _{L^2} ^2 - \Vert D^{m-1} \uu(T) \Vert _{L^2} ^2 \\ & - 4 \chi ^2 \mu ^{-1} \int _t ^T \Vert D^{m-1} \error (\tau) \Vert _{L^2} ^2 \, d \tau,
\end{align*}
for $T > t > t_{\ast}$.  Taking a $t$ larger than the time for which the induction step is valid,  and also so that \eqref{eqn:zeta-m} in Theorem \ref{monotonicity-zeta} and \eqref{eqn:upper-bound-du} in Theorem \ref{decay-du-dw-upper} hold,  we have  
 
\begin{align*}
4 \mu \, T   \Vert D^m \uu (\tau) \Vert _{L^2} ^2 &  \geq C \mu ^{-(m-1)} \left(t^{- (2 \alpha + m - 1)} - \lambda_m ^2 \, T^{\, -(2 \alpha + m - 1)} \right) \\ & - 4 \chi ^2 \mu ^{-1} \int _t ^{\infty} \Vert D^{m-1} \error (\tau) \Vert _{L^2} ^2 \, d \tau - 16 \chi ^2 \, T \int _t ^{\infty} \Vert D^{m} \error (\tau) \Vert _{L^2} ^2 \, d \tau,
\end{align*}
where $\lambda_m = \frac{C_{m-1}}{\widetilde{C}_{m-1}} \geq 1$,  where $C_{m-1}$ is the constant in \eqref{eqn:upper-bound-du} in Theorem \ref{decay-du-dw-upper} and $\widetilde{C}_{m-1}$ is the constant in the induction step.  Taking $T = \left(2 \lambda_m \right) ^{\frac{1}{\alpha}}  t$ and using \eqref{eqn:first-estimate-dmerror} from Theorem \ref{estimate-error-first}, we obtain the result. $\Box$

\begin{Remark}
As a consequence of \eqref{eqn:upper-bound-du},  \eqref{eqn:first-estimate-dmerror},  \eqref{eqn:equation-u-rewritten} and \eqref{eqn:lower-bound-du}, we can show that $\Vert D^m \uu (t) \Vert _{L^2}$ is monotonically decreasing in $(\zeta _m ^{'},   \infty)$, for some large enough $\zeta _m ^{'}$.  This,  improves \eqref{eqn:zeta-m} in Theorem \ref{monotonicity-zeta}.  Moreover,  from \eqref{eqn:error},  \eqref{eqn:equation-u-rewritten} and \eqref{lower-bound-du}, we have that for $m \geq 0$

\begin{displaymath}
\Vert D ^m \ww (t) \Vert _{L^2} \geq C \mu ^{- \frac{m+1}{2}} \, t^{- \left(\alpha + \frac{m+1}{2}  \right)},
\end{displaymath}
for large enough $t$.  Also, when $\mu = \nu$,  from \eqref{eqn:upper-bound-du},  \eqref{eqn:upper-bound-dw},   \eqref{eqn:upper-bound-error}, \eqref{eqn:lower-bound-dw} and  \eqref{eqn:equation-w-rewritten},   we deduce that $\Vert D ^m \ww (t) \Vert _{L^2}$ is monotonically decreasing in $(\zeta _m ^{''},   \infty)$, for some large enough $\zeta _m ^{''}$. 
\end{Remark}

\bibliographystyle{plain}
\bibliography{GMNPZ}{}

\end{document}